\documentclass[11pt]{article}
\input{amssymb.sty}
\usepackage[pctex32]{graphics}
\input{psfig.sty}
\begin{document}
\bibliographystyle{unsrt}
\newtheorem{thm}{{\sc Theorem}}[section]
\newtheorem{lma}{{\sc Lemma}}[section]
\newtheorem{prp}{{\sc Proposition}}[section]
\newtheorem{cor}{{\sc Corollary}}

\newfont{\msbm}{msbm10 scaled\magstep1}
\newfont{\eusm}{eusm10 scaled\magstep1}
\newfont{\goth}{eufm10 scaled\magstep1}
\newfont{\lasym}{lasyb10 scaled\magstep1}
\newfont{\lasys}{lasy7 scaled\magstep1}

\def\diam{{\mbox{\lasym 3}}}
\def\diams{{\mbox{\lasys 3}}}

\def\bea*{\begin{eqnarray*}}
\def\eea*{\end{eqnarray*}}
\def\ba{\begin{array}}
\def\ea{\end{array}}
\count1=1
\def\be{\ifnum \count1=0 $$ \else \begin{equation}\fi}
\def\ee{\ifnum\count1=0 $$ \else \end{equation}\fi}
\def\ele(#1){\ifnum\count1=0 \eqno({\bf #1}) $$ \else \label{#1}\end{equation}\fi}
\def\req(#1){\ifnum\count1=0 {\bf #1}\else \ref{#1}\fi}
\def\bea(#1){\ifnum \count1=0   $$ \begin{array}{#1}
\else \begin{equation} \begin{array}{#1} \fi}
\def\eea{\ifnum \count1=0 \end{array} $$
\else  \end{array}\end{equation}\fi}
\def\elea(#1){\ifnum \count1=0 \end{array}\label{#1}\eqno({\bf #1}) $$
\else\end{array}\label{#1}\end{equation}\fi}
\def\cit(#1){
\ifnum\count1=0 {\bf #1} \cite{#1} \else 
\cite{#1}\fi}
\def\bibit(#1){\ifnum\count1=0 \bibitem{#1} [#1    ] \else \bibitem{#1}\fi}
\def\ds{\displaystyle}
\def\hb{\hfill\break}
\def\comment#1{\hb {***** {\em #1} *****}\hb }

\newcommand{\TZ}{\hbox{\bf T}}
\newcommand{\MZ}{\hbox{\bf M}}
\newcommand{\NZ}{\hbox{\bf N}}
\def\lx{{\prec_{\rm L}}}
\def\lt{{\rm lt}}
\def\vv{{\bf v}}
\def\vr{{\bf r}}
\def\vp{{\bf p}}
\def\spec{{\rm spec}}
\def\Hom{{\rm Hom}}
\def\eM{{{J_0^\bot}}}
\def\Io{{I(\widehat{o})}}
\def\wt{{\rm wt}}
\def\RT{{\cal R}}
\def\hl{{\rm  Hilb} }
\def\NZ{{\mathbb N}}
\def\ZZ{{\mathbb Z}}
\def\RZ{{\mathbb R}}
\def\CZ{{\mathbb C}}
\def\PZ{{\mathbb P}}
\def\QZ{{\mathbb Q}}
\def\HB{{\hl^{\vert G\vert }(\CZ^n)^G}}
\def\HBA{{\hl^{\vert A_1(5)\vert }(\CZ^5)^{A_1(5)}}}
\def\HBB{{\hl^{\vert A_1(4)\vert }(\CZ^4)^{A_1(4)}}}
\def\HS{{\hl^G(\CZ^n)}}
\def\HSA{{\hl^{A_1(5)}(\CZ^5)}}
\def\HSB{{\hl^{A_1(4)}(\CZ^4)}}
\newcommand{\HZ}{\hbox{\bf H}}
\newcommand{\EZ}{\hbox{\bf E}}
\newcommand{\GZ}{\,\hbox{\bf G}}
\newtheorem{pf}{Proof}
\renewcommand{\thepf}{}
\vbox{\vspace{38mm}}

\def\lra{{\leftrightarrow}}
\begin{center}
{\LARGE \bf  Crepant Resolutions of $\CZ^n/A_1(n)$ and Flops of $n$-Folds for $n=4,5$ }\\[5mm]
Li Chiang 
\\{\it Institute of Mathematics \\ Academia Sinica \\ 
Taipei , Taiwan \\ }
(e-mail:chiangl@gate.sinica.edu.tw)
\\[5mm]
Shi-shyr Roan
\\{\it Institute of Mathematics \\ Academia
Sinica \\  Taipei , Taiwan \\ } 
(e-mail:maroan@ccvax.sinica.edu.tw) \\[5mm]
\end{center}

\begin{abstract} 
In this article, we determine the explicit toric variety structure of $\hl^{A_1(n)}(\CZ^n)$ for $n=4,5$, where $A_1(n)$ is the special diagonal group of all order $2$ elements. Through the toric data of $\hl^{A_1(n)}(\CZ^n)$, we obtain certain toric crepant resolutions of $\CZ^n/A_1(n)$, and the different crepant resolutions are connected by flops of $n$-folds for $n=4,5$.
\par \vspace{5mm} \noindent
2000 MSC Primary: 14M25,14J17, 20C33;  Secondary:13P10

\end{abstract}
\vfill
\eject




\begin{section}{Introduction}
This article is intended to discuss the construction of flops of high-dimensional manifolds through the geometry of orbifolds. By an orbifold, we mean an algebraic (or analytic) variety with at most quotient singularities. This geometrical object was
introduced in the 1950s by
I. Satake, which he called a V-manifold, 
on his course of study of the
Siegel's modular forms
\cite{S}.
Special emphasis in this paper will be on the resolutions of orbifolds with 
the Gorenstein quotient singularities, i.e., the local structure around a singular point is given by 
\begin{equation}
\CZ^n/ G \ , \ \ \ \  G : {\rm a \ (nontrivial) \ 
subgroup \ of \ } {\rm SL}_n ( \CZ ) \ , \ n \geq 2 \ .
\label{localstr}
\end{equation}
In the mid-80s,
string theorists studied a special type of 
physical model where strings propagate on an orbifold $ M/G$ where 
$M$ is a manifold quotiented by a finite group  
$G$.  The required  
string vacuum  is described by the following orbifold Euler characteristic  \cite{DHVW}:
$$
\chi(M, G) = \frac{1}{|G|} \sum_{
g, h \in G, \ gh=hg }\chi(M^{g, h}) \  
$$
where $M^{g, h}$ is the simultaneous fixed-point set of $g, h$. 
An equivalent expression of the above formula is 
given by 
$$
\chi(M, G) = \sum_{ [g] }\chi(M^g/C(g)) , 
$$
where $[g]$ is the conjugacy class of $g$ with the centralizer $C(g)$. 
In the local situation (\ref{localstr}), one obtains 
$$
\chi(\CZ^n , G) = | {\rm Irr}(G) | \ , \ \ 
{\rm Irr}(G) : = \{ {\rm irruducible \ representations \ of \ } G \} \ .
$$
For the consistency of physical theory,  the Gorenstein orbifold $\CZ^n/G$ demands a resolution $\widehat{\CZ^n/G}$ with the
trivial canonical bundle and the Euler number
$\chi(\widehat{\CZ^n/G}) =\chi(\CZ^n, G)$. For $n=2$, all finite subgroups $G$ of ${\rm SL}_2(\CZ)$ were classified into $A$-$D$-$E$ series by F. Klein, where the orbifolds $\CZ^2/G$ are always of the hypersurface singularity. The minimal resolution of $\CZ^2/G$ in the surface theory  provides the required crepant resolution $\widehat{\CZ^2/G}$ \cite{HH}.
For $n=3$, the desired crepant resolution $\widehat{\CZ^3/G}$ was found for all $G$ in the mid-nineties \cite{I1, Ma, R94, Rtop} through some quantitative, but constructive, methods, heavily depending on the classical  Miller-Blichfelt-Dickson classification of finite subgroups of ${\rm SL}_3(\CZ)$ \cite{MBD}, and the Klein's invariant theory of the simple groups $I_{60}$ (Icosahedral group) and $H_{168}$ (the simple group of  order 168) \cite{R94}. The non-uniqueness of such crepant resolution is expected due to the flop-relation in $3$-folds, but all with the required Euler number. Thereafter, the qualitative understanding and a possible canonical one of these crepant resolutions naturally arise as problems in algebraic geometry. For this purpose, the development has resulted in the concept of Hilbert scheme of $G$-orbits (or $G$-Hilbert scheme), $\HS$, associated to $\CZ^n/G$, which was proposed in \cite{IN1} by engaging the finite-group representation theory into geometrical problems of orbifolds. The structure of $G$-Hilbert scheme has now provided a positive answer for the problem addressed to the crepant resolution with the orbifold Euler number for the dimension $n=2, 3$ \cite{BKR, INj, Na}.     

For $n \geq 4$, very few results have been known on the crepant resolution problem of $\CZ^n/G$, and the structure of $\hl^G(\CZ^n)$, partly due to the lack of common conclusion for all $G$, as $\CZ^n/G$ has no crepant resolution at all for a large class of groups $G$, even in the abelian group case. A suitable class of groups $G$, which could be of interest to study, are those with Gorenstein orbifold $\CZ^n/G$ of the {\it hypersurface type}, and one can regard it as one way to the high-dimensional generalization of Klein surface singularities. For the abelian group case, the group $G$ is given by 
$$A_r(n):=\{g\in\mbox{SL}_n(\CZ)\vert  g^{r+1}=1,\ g \mbox{ is diagonal}\} \  \ , \ \ r \geq 1 \ .
$$
In \cite{CR1}, we have studied the cases $A_r(4)$ for all $r$, and obtained the crepant resolutions $\widehat{\CZ^4/A_r(4)}$ through the detailed structure of $\hl^{A_r(4)}(\CZ^4)$.  
While in this work, we shall discuss the cases $A_1(n)$ for $n=4,5$, in order to explore the flop-relation of $n$-folds.
It is known that there do exist some crepant resolutions of $\CZ^n/A_1(n)$ for all $n$ in \cite{R94}, where the result was obtained by blowing-up techniques through the hypersurface equation of $\CZ^n/A_1(n)$. But in this paper, we shall construct another type of crepant resolutions of $\CZ^n/A_1(n)$, which are toric varieties, through the toric structure of $\hl^{A_1(n)}(\CZ^n)$. One of the main goals of this approach is to identify the proper ``flop'' concept of $n$-folds for $n=4,5$. This method is consistent with the flops of $3$-folds, considered as crepant resolutions of the $3$-dimensional isolated singularity with the equation, $XY=ZW$, a construction known since the seventies in the study of degeneration of K3 surfaces \cite{R0}. However, for $n=4,5$, the singularities, which possess crepant resolutions with the flop relation discovered in this article, has a complicated expression of algebraic equations (see formula (\ref{4sing}), (\ref{5sing}) of the paper). In particular, it is not of the type of complete intersection.

This paper is organized as follows: In \S2, we give a brief survey on the construction of  crepant resolution $\widehat{\CZ^n/G}$ for $n=2, 3$.   In \S3, we review some basic facts in $G$-Hilbert schemes and toric geometry for later discussions. In \S4, we derive the toric variety structure of $\HSB$, which is smooth but not crepant. By blowing down the canonical divisor $(\PZ^1)^3$ of $\HSB$ to $(\PZ^1)^2$ in different factors, we obtain three toric crepant resolutions of $\CZ^4/A_1(4)$. By this process, the flop of $4$-folds is naturally revealed between these three crepant resolutions. In \S5, we derive the toric variety structure of $\HSA$, which is singular and non-crepant. By analyzing the toric structure of $\HSA$, we construct twelve crepant resolutions of $\CZ^5/A_1(5)$, all dominated by $\HSA$. The connection between these crepant resolutions gives rise to the ``flop'' relation of $5$-folds.  
\end{section}

\begin{section}{Kleinian Singularities and Crepant Resolutions of Gorenstein Three-Dimensional Orbifolds } 
For a finite non-trivial subgroup $G$ of ${\rm SL}_2 ( \CZ )$, $\CZ^2 / G$ is a surface with 
an isolated singularity at the origin. These singularities were classified by Klein in 1872 ( or so ) in his work on the invariant theory of regular solids in 
$\RZ^3$, and they are described in the following table:

{\center
\begin{tabular}{| c |c |c| c| c| }
\hline 
{\rm\small Type }& {\rm\small Group} $G$ & $\vert G \vert$ &{\rm\small degree of invariants} &{\rm\small relation \ of \ invariants} \\ 
\hline
$A_r$ & {\small\rm Abelian } & $r+1$ & $2, r+1, r+1$ & $X^{r+1} + YZ$ \\
$D_r$ & {\small\rm Binary dihedral} & $4(r-2)$ & $4,2(r-2),2(r-1)$ & $X^{r-1} + XY^2 + Z^2$ \\
$E_6$ & {\small\rm Binary tetrahedral} & $24$ & $6, 8, 12$ & $X^3 + Y^4 + Z^2$ \\
$E_7$ & {\small\rm Binary octahedral} & $48$ & $8, 12, 18$ & $X^3Y + Y^3 + Z^2$ \\
$E_8$ & {\small\rm Binary icosahedral} & $120$ & $12, 20, 30$ & $X^5 + Y^3 + Z^2$ \\
\hline 
\end{tabular}}

\noindent
For $A$-type group $G$, it is a cyclic 
subgroup 
of ${\rm SL}_2(\CZ)$ generated by the diagonal element with the entries 
$e^{\frac{2\pi i }{r+1}} , 
e^{\frac{2\pi i r}{r+1}}$. By the continued fraction, 
$$
\frac{r+1}{r} = 2 - \frac{1}{2 - \frac{1}{\ddots}} \ \ ,
$$ 
the minimal resolution 
$\widehat{\CZ^2/G}$ has the exceptional set consisting of a tree of $r$-rational curves with 
self-intersections $-2$, which implies   
the trivial canonical bundle of $\widehat{\CZ^2/G}$ with the Euler number $r+1$. For other groups $G$ in 
${\rm SL}_2(\CZ)$, the minus identity map always belongs to the center of $G$. 
One can construct the minimal resolution by first blowing up 
$\CZ^2$ at the origin, then  analyzing the fixed-points of $G$ on this 
blowing-up space, and then reducing to cyclic quotient singularities. 
The following remarkable 
characterization for Kleinian  singularities by Brieskorn is known:
\par \vspace{0.2mm} \noindent
\begin{thm} \label{th:B}
 \cite{B} For an isolated normal surface 
singularity, the following are 
equivalent: 
(1) it is Kleinian , (2) it is rational with multiplicity 2, (3) it is rational with embedding dimension 3, and (4) the minimal resolution has the exceptional set 
by a configuration of rational (-2)-curves with an intersection matrix 
minus that of the Partan matrix 
of the root system of type $A_r, D_r, E_6, E_7$ or $E_8$. The correspondence 
between (1) and (4) is indicated in the above table.
\end{thm} 
On the other hand, J. McKay observed an interesting
 connection between the representation theory of Kleinian groups and the affine 
Dynkin diagrams of the $A$-$D$-$E$ root system. Denote $c : G \longrightarrow {\rm GL}_2 ( \CZ )
$  the `canonical' representation of $G$. For $\rho, \rho' \in {\rm Irr}(G)$, let 
$m (\rho, \rho')$ be the multiplicity of $\rho'$ in the representation 
$\rho \otimes c$:
\[ 
\rho \otimes c = \sum_{\rho' \in {\rm Irr}(G)} m (\rho, \rho') \rho' \ .
\]    
By the self-dual property of $c$,  $m (\rho, \rho') = m (\rho', \rho)$ for all 
$\rho, \rho'$. One forms a diagram $\widetilde{\Gamma} ( G )$ with 
elements in Irr(G) as vertices, and two vertices $\rho, \rho'$ 
connected by $m ( \rho, \rho')$ edges. Then, one has the following 
result:
\par \vspace{0.2mm} \noindent 
\begin{thm} \label{th:M}
 \cite{Mc} \ 
$\widetilde{\Gamma} ( G ) $ is isomorphic to the affine 
Dynkin diagram of the corresponding root system.  
 \end{thm} 
By a suitably chosen isomorphism between $\widetilde{\Gamma} ( G ) $ and the affine 
Dynkin diagram, the trivial representation of $G$ corresponds to the negative 
longest root of the affine diagram. So, the diagram $\Gamma ( G )$ 
obtained by deleting 
the trivial representation is isomorphic to the ordinary Dynkin diagram of 
its corresponding system. With the results in Theorems \ref{th:B} and \ref{th:M}, the rational 
(-2)-curves in resolution configuration implies 
the trivial canonical bundle of $\widehat{\CZ^2/G}$, and the Euler number equality
follows from McKay's correspondence \cite{HH}. 

For $n=3$, all finite subgroups $G$ of ${\rm SL}_3 (\CZ )$
were classified by Miller et al. in 1916  
\cite{MBD}. The list consists of 12 types of such groups, 
among which two well-known simple groups, $I_{60}$  
and $H_{168}$  appear. For the problem of  the
crepant resolutions of $\widehat{\CZ^3/G}$, due to the long list of 
the  Miller-Blichfelt-Dickson classification, it is still hard to 
study the problem by the method of case by case approach. Hence,
certain induction procedures are introduced, as the one of blowing up at the origin in 
the 2-dimensional 
case,  to enable us to 
 work only on a few cases, consisting of abelian groups, certain solvable 
groups and simple groups $I_{60}, H_{168}$ \cite{R94}. For 
abelian groups, the toric 
geometry \cite{KKMS} provides the effective computational methods to obtain the solutions \cite{MOP, RY, R89}. Using these techniques, one can work on 
an explicitly given solvable subgroup of ${\rm SL}_3(\CZ)$ and obtain  
the crepant resolutions by using a certain `canonical' procedure, regardless the resulting structure of the resolution might be complicated. For the  
simple groups $I_{60}, H_{168}$ in ${\rm SL}_3(\CZ)$, one needs some different 
techniques, using the invariant theory of their 3-dimensional 
representations in Klein's work \cite{Klein}:
$$
\begin{array}{| c | c | c |}
\hline 
{\rm G} & {\rm degree \ of \ invariants} & {\rm relation \ of \ invariants}\\
\hline 
I_{60} & 2, 6, 10, 15  &-144W^2 + Z^3 - 1728 Y^5 + 720 XY^3Z \\
& & - 80 X^2YZ^2 + 64X^3(5Y^2-XZ)^2  \\
\hline
H_{168} &  4, 6, 14, 21 &-W^2 + Z^3 + 1728 Y^7 + 1008 XY^4Z \\ 
& & - 88 X^2YZ^2 -60032X^3Y^5+ 1088X^4Y^2Z   \\
& & + 22016X^6Y^3 - 256 X^7Z - 2048 X^9Y \\
\hline
\end{array}
$$
Note that the relation
of invariants in the above table define a hypersurface
of $\CZ^3/G$ in $\CZ^4$. One can obtains the required crepant resolution by making full use of the relations of invariants of these simple groups $G$ \cite{Ma, R94}. With all these efforts, one 
arrives at a solution of the crepant resolutions 
$\widehat{\CZ^3/G}$ for a finite subgroup $G$ of ${\rm SL}_3(\CZ)$  \cite{Rtop}. However, all of those 
works were carried out using an explicit (somewhat unpleasant) 
computational method, which is cumbersome and not very illuminating. An approach of a qualitative nature to the crepant resolution $\widehat{\CZ^3/G}$ is provided by the method of $G$-Hilbert scheme, which will be briefly described in the next section.

\end{section}

\begin{section}{$G$-Hilbert scheme, Abelian Orbifolds and Toric Geometry} 
From now on for the rest of this paper, $\{e^1,..,e^n\}$ will always denote the standard basis of $\CZ^n$, $\{Z_1,..,Z_n\}$ the corresponding dual basis, and $\CZ[Z]:=\CZ[Z_1,..,Z_n]$ the polynomial ring of $Z_i$'s. 

For a finite subgroup $G$ of $\mbox{GL}_n(\CZ)$, we denote $S_G=\CZ^n/G$ with the canonical projection $\pi_G:\CZ^n\rightarrow S_G$,
and $o=\pi_G(\vec{0})$. In this article, a variety $X$ is said to be birational over $S_G$ if there exists a proper birational epimorphism $\sigma_X:X\rightarrow S_G$. In this case, we have the commutative diagram: 
$$\matrix{{X \times_{S_G} \CZ^n}  & \longrightarrow & \CZ^n \cr
\ \ \downarrow \pi & & \downarrow \pi_G  \cr
X  &\stackrel{\sigma_X}{\longrightarrow} & S_G  \
.
}
$$
Denote ${\cal F}_X\left(=\pi_* {\cal O}_{X \times_{S_G}\CZ^n}\right)$
the coherent ${\cal O}_X$-sheaf over $X$ obtained by 
the push-forward of the structure sheaf of $X \times_{S_G} \CZ^n$. 
The geometrical fiber of ${\cal F}_X$ over an element $y$ of $X$ is given by
${\cal F}_{X, y} = k(y) \bigotimes_{{\cal O}_X} {\cal F}_X$, which is isomorphic to $\CZ[Z]/I(y)$ for some $G$-invariant ideal $I(y)$ in $\CZ[Z]$. 
The $G$-Hilbert scheme $\hl^G(\CZ^n)$ with the birational morphism,
$\sigma_\hl  : \HS \longrightarrow  S_G$, 
is the minimal object in the category of varieties $X$ birational over $S_G$ such that 
${\cal F}_X$ is a vector bundle over $X$ \cite{IN1} (or see \cite{CR1} and references therein). Furthermore, for each $X$ in the category, there exists a unique epimorphism $\lambda_X$ from $X$ to $\HS$. Indeed, $\lambda_X(y)$ is represented by the ideal $I(y)$. In the case when  $C^n/G$ in (\ref{localstr}), for $n=2$, $\hl^G(\CZ^2)$ is  the minimal resolution of $\CZ^2/G$, and for $n=3$, $\hl^G(\CZ^3)$ is a toric crepant resolution of $\CZ^3/G$ for an abelian group $G$ \cite{INj, Na}. For a general finite subgroup $G$ of ${\rm SL}_3(\CZ)$, $\hl^G(\CZ^3)$ is indeed a  crepant resolution $\widehat{\CZ^3/G}$ has been justified in \cite {BKR} by methods in homological algebra, an argument bypassing the geometry of $G$-Hilbert scheme (of which explicit structures are important on its own from the aspect of applications to certain relevant physical problems).

When $G$ is an abelian group, one can apply toric geometry to study problems on resolutions of $S_G$. We now give a brief review of some facts in toric geometry for later use in this article (for the details, see \cite{KKMS}). For the rest of this paper, we shall always assume $G$ to be a finite abelian group in ${\rm SL}_n(\CZ)$, and identify $G$ as a subgroup of the diagonal group $T_0:={\CZ^*}^n\subset\mbox{GL}_n(\CZ)$.
Regarding $\CZ^n$ as the partial compactification of $T_0$, we
define the $n$-torus $T$ with the embedding in the $T$-space $S_G$ by 
$$T := T_0/G \ ,\ \ \ T \subset S_G \ .$$
Denote $N = {\rm Hom}(\CZ^*, T )$ (resp. $N_0 = {\rm Hom}(\CZ^*, T_0 )$) the lattice of one-parameter subgroups of $T$ (resp. $T_0$), $M$ (resp. $M_0$) the dual lattice of $N$ (resp. $N_0$). We have $N \supset N_0$ and $M  \subset M_0$. Through the map ${\rm exp}: {\RZ}^n \longrightarrow T_0 $, ${\rm exp}(\sum_i x_i e^i) := \sum_i e^{2\pi\sqrt{-1}x_i}e^i$, one can identify $N_0$ and $N$ with the following lattices in $\RZ^n$,
$$
N_0  = \ZZ^n (= {\rm exp}^{-1}(1)) \ , \ \ 
N = {\rm exp}^{-1}( G ) \ .
$$
In this way, $M_0$ can be identified with the monomial group of variables $Z_1,..,Z_n$, and $M$ corresponds to the subgroup of $G$-invariant monomials in $M_0$.

By \cite{KKMS}, a $T$-space birational over $S_G$ is described by a fan $\Sigma = \{ \sigma_l \ | \ l \in I\}$ with 
$C_0:=\sum_i\RZ_{\ge 0}e^i$
as its support, i.e., a rational polyhedral cone decomposition of $C_0$. In this situation, an equivalent description is given by the polytope decomposition 
$\Lambda  = \{\Delta_l \ | \ l \in I \}$ of the simplex
$$
\Delta:= \{ \sum_ix_ie^i \in C_0 | \sum_ix_i=1 \},
$$
where $\Delta_l:= \sigma_l \bigcap \Delta$ with the vertices of $\Delta_l$ in $\Delta\cap \QZ^n$. For $\sigma_l=\{{\bf 0}\}$, we have $\Delta_l=\emptyset$. We shall call $\Lambda$ a rational polytope decomposition of $\Delta$, and denote $X_\Lambda$ the toric variety corresponding to $\Lambda$.

For a rational polytope decomposition $\Lambda$ of
$\Delta$, we define $\Lambda(i):= \{ \Delta_l \in \Lambda \ | \ {\rm dim}(\Delta_l) = i \}$ 
for $ -1 \leq i\leq n-1$, (here ${\rm dim}( \emptyset ):= -1$). Then, for each $\Delta_l \in \Lambda(i)$, there associates a $T$-orbit of dimension $n-1-i$, denoted by ${\rm orb}(\Delta_l)$. For $\Delta_l\in \Lambda(n-1)$, i.e., $\Delta_l$ is an $(n-1)$-dimensional polytope, ${\rm orb}(\Delta_l)$ consists of only one point, which will be denoted by $x_{\Delta_l}$. For a vertex $v\in \Lambda(0)$, the closure of ${\rm orb}(v)$ represents a toric divisor, which will be denoted by $D_v$. The canonical sheaf of $X_\Lambda$ is given by the following expression of toric divisors,
\begin{equation}
\omega_{X_\Lambda} = {\cal O}_{X_\Lambda}(\sum_{v \in \Lambda(0) } (m_v-1)D_v ) \ , \ 
\label{cb}
\end{equation}
where $m_v$ is the smallest positive integer with $m_v v\in N$, i.e., a primitive element of $N$. 

A polytope decomposition $\Lambda$ of $\Delta$ is called integral if all vertices of $\Lambda$ are in $N$. By (\ref{cb}), the trivial canonical sheaf is described by the integral condition of $\Lambda$. The non-singular criterion of $X_{\Lambda}$ is given by the simplicial decomposition of $\Lambda$ with the multiplicity-one property. For a polytope $\RT\in \Lambda(i)$, we denote 
$$X_{\RT}=\spec(\CZ[M \cap \check{\sigma_{\RT}}])$$
where $\check{\sigma_{\RT}}$ is the dual of the cone ${\sigma_{\RT}}:=\{r v\vert r\ge 0,\ v\in \RT \}.$
Then, $\{X_\RT\}_{\RT\in \Lambda(n-1)}$ forms an open cover of $X_\Lambda$. 

Now we consider the case $G=A_1(n)$. We shall always denote the following elements in $M(=M_0^{G})$:
\begin{equation}
T_i=Z_i^2, \ \ \ \ \ U_i={{Z_i}\over{(Z_1..\hat {Z_j}..Z_n)}},\ \ \ \ \   V_{ij}={{Z_iZ_j}\over{\prod_{\alpha\ne i,j;\ 1\le\alpha\le n}Z_\alpha}}\ .
\label{eq:TUV}
\end{equation}
The set of $N$-integral elements in $\Delta$ is given by 
$$\Delta\cap N=\{e^i\vert 1\le i\le n\}\cup \{v^{ij}\vert 1\le i< j\le n\},  \ \ \ \ \ \ \ v^{ij}:={1 \over 2}(e^i+e^j).$$
Other than the simplex $\Delta$ itself, there is only one integral polytope decomposition of $\Delta$, denoted by $\Xi$, which is invariant under all permutations of the coordinates (i.e., $\tau:e^i\mapsto e^{\tau(i)}$ for $\tau\in {\goth S}_n$, the symmetric group of degree $n$). Then, $\Xi(n-1)$ consists of the following $n+1$ polytopes: 
\begin{equation}
\Diamond:=\langle v^{ik}\vert 1\le i< k\le n\rangle,\ \ \ \ \ \ \ \ \ \Delta_j:=\langle e^j, v^{ij}\vert i\ne j\rangle \ \ \ \ (j=1,..,n).
\label{def:DD}
\end{equation}
The toric variety $X_{\Xi}$ is smooth on the affine space centered at $x_{\Delta_j}$ with the local coordinate $U_j$, $T_i$ $(i\ne j)$. Therefore, the singular set of $X_{\Xi}$ lies in the affine open set $X_{\Diamond}$ centered at $x_{\Diamond}$. The structure near $x_{\Diamond}$ is determined by $\check{\sigma_\Diamond} \cap M$, which is generated by
$T_i,\ \ U_i^{-1}.$
For $n=3$, $U_1^{-1}$, $U_2^{-1}$, $U_3^{-1}$ form an integral basis for $M$, which implies the smoothness of $X_\Xi$. For $n\ge 4$, $X_{\Xi}$ is singular at $x_\Diamond$, and $T_i$, $U_i^{-1}$, $1\le i\le n$, form the minimal generating set of $\check{\sigma_\Diamond} \cap M$. For $n=4$, the singular structure of $X_{\Xi}$ near $x_\Diamond$ is described by the $4$-dimensional affine variety in $\CZ^8$ with the equations:
\begin{eqnarray}
t_i\overline{u}_i= t_j\overline{u}_j \ , \ \ \  t_it_j = \overline{u}_{i^\prime}\overline{u}_{j^\prime} \ , \ \ \ \ \ \ \  (t_i, \overline{u}_i)_{1\leq i \leq 4} \in \CZ^8
\label{4sing}
\end{eqnarray}
where $i \neq j$ and $\{ i^\prime, j^\prime \} $ is the
complementary pair of $\{i, j\}$. For $n=5$, the singular structure of $X_\Xi$ near $x_\Diamond$ is defined by the following relations: 
\begin{equation}
t_i\overline{u}_i= t_j \overline{u}_j ,
\ \
\ \overline{u}_i\overline{u}_j= t_kt_lt_m  \ ,  \ \ \ \ \ \ \  (t_i, \overline{u}_i)_{1\leq i \leq 5} \in \CZ^{10} 
\label{5sing}
\end{equation}
where $1\le i< j\le 5$, and $\{k,l,m\}$ is the complement of $\{i,j\}$. 

For later use, we recall some terminology of the Gr\"{o}bner basis \cite{CLO}. Let $\prec_{\rm L}$ be a lexicographic order on $\CZ[Z]$ and ${\bf w}\in \ZZ_{\ge 0}^n$. The weight order $\prec$ determined by the weight ${\bf w}$ is the monomial ordering on $\CZ[Z]$ with the following properties:
 
\begin{tabular}{c l}
(1)& If $({\bf u}_1-{\bf u}_2)\cdot {\bf w}<0$, then $Z^{{\bf u}_1}\prec Z^{{\bf u}_2}$.\\
(2)& If $({\bf u}_1-{\bf u}_2)\cdot {\bf w}=0$ and $Z^{{\bf u}_1}\prec_{\rm L}Z^{{\bf u}_2}$, then $Z^{{\bf u}_1}\prec Z^{{\bf u}_2}$.
\end{tabular}

\noindent For a monomial ideal $I\subset \CZ[Z]$, we shall denote  $I^\bot$ the set of monic monomials outside $I$. For a monomial ordering $\prec$ and an ideal $J$, $\lt_\prec(J)$ will denote the ideal generated by the leading terms (with respect to $\prec$) of all elements in $J$.

\end{section}

\begin{section}{$\HSB$ and Crepant Resolutions of $\CZ^4/A_1(4)$}
In this section, we consider the case $n=4$. Throughout this section, the indices $j,k,l,m$ will always denote a permutation of $(1,2,3,4)$. For $n=4$, the polytope $\Diamond$ in (\ref{def:DD}) is a regular $3$-dimensional octahedron contained in the standard simplex $\Delta$ in $\RZ^4$, with a cube as its dual polygon. We label the facets of the octahedron $\Diamond$ by $F_j, F_j^\prime$ where 
$$
F_j = \Diamond \cap \bigtriangleup_j \ , \ \ \
F_j^\prime = 
\{ \sum_{i=1}^4 x_ie^i \in \Diamond \ | \ x_j =0
\} \ .
$$
Then, the dual of $F_j$, $F_j^\prime$ are vertices of the cube, denoted by $\alpha_j$, $\alpha_j^\prime$ (see [Fig. 1]).
\vskip 10pt
\begin{tabular}{c}
$\mbox{\psfig{figure=octeface.ps}}$ \hskip 70pt
$\mbox{\psfig{figure=cubic.ps}}$\\
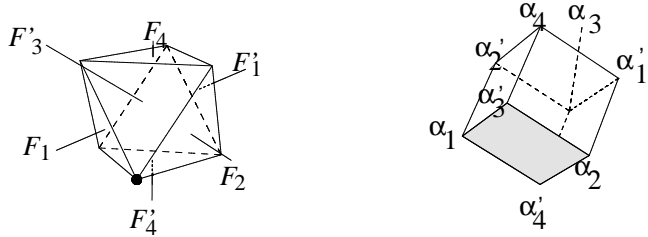
{[Fig. 1] Dual pair of octahedron and cube: faces $F_j$, $F_j^\prime$ of octahedron }\\
{are dual to vertices $\alpha_j$, $\alpha_j^\prime$ of cube.}
\end{tabular}

\begin{tabular}{c}
$\mbox{\psfig{figure=tatrhilb.ps}}$\\
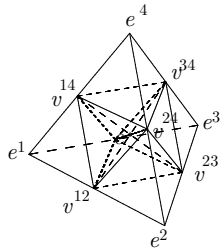
{[Fig. 2] The rational simplicial decomposition $\Xi^*$ of $\Delta$ for $n=4$.}
\end{tabular}
\vskip 10pt 

\noindent Consider the rational simplicial decomposition $\Xi^*$ of $\Delta$ obtained by $\Xi$ and adding the center 
$c := \frac{1}{4} \sum_{j=1}^4 e^j$
as a vertex of $\Xi^*$ with the barycentric decomposition
of $\Diamond$ in $\Xi$, (see [Fig. 2]). Then, $\Xi^*(3)$ consists of the polytopes $\Delta_j$, $C_j$ and $C^\prime_j$ for $j=1,..,4$, where
$$C_j:=\langle c, v^{jk}, v^{jl}, v^{jm}  \rangle,\ \ \ 
C^\prime_j:=\langle c, v^{kl}, v^{lm}, v^{mk} \rangle.$$ 
Note that $c\not\in N$ and $2c\in N$.
\begin{thm}\label{th:A1(4)}
For $G= A_1(4)$, we have 
$\HS  \simeq 
X_{\Xi^*} \ $, 
which is non-singular with the canonical bundle 
$\omega = {\cal O}_{X_{\Xi^*}} (E)$, where $E$ is
an irreducible divisor isomorphic to the triple
product of $\PZ^1$, 
$E =\left( \PZ^1 \right)^3$. Furthermore, for $\{j,k,l \} = \{1,2,3\}$ and the $(k,l)$-th factor projection, $p_j:E\longrightarrow\left(\PZ^1\right)^2$, the restriction of normal bundle of $E$ on each fiber $(\simeq
\PZ^1)$ of $p_j$  is the $(-1)$-hyperplane bundle of $\PZ^1$.
\end{thm}
\noindent{\it Proof}: 
The smoothness of the affine spaces $X_{\RT}$, $\RT=\Delta_j, C_j, C_j^\prime$, follows from the integral and multiplicity-one properties of $\sigma_\RT$ with respect to $N$. With (\ref{eq:TUV}) for $n=4$, we have $V_{jk}=(Z_jZ_k)/(Z_lZ_m)$, $U_j=Z_j/(Z_kZ_lZ_m)$. By using the inverse matrices of $\left(e^j,v^{jk}, v^{jl}, v^{jm}\right)$, $\left(2c, v^{jk}, v^{jl}, v^{jm}\right)$ and $\left(2c, v^{kl}, v^{lm}, v^{mk}\right)$,
one obtains the generators of $\check{\sigma_{\RT}}\cap M$ for $\RT=\Delta_j, C_j$, and $C_j^\prime$ as follows:
$$\matrix{
{\Delta_j:\ \ \ \ \ \ }&{ U_j, }&{T_k, }&{T_l, }&{T_m ;}\cr
{C_j:\ \ \ \ \ \ }&{U_j^{-1}, }&{V_{jk}, }&{V_{jl}, }&{ V_{jm} ;}\cr
{C_j^\prime:\ \ \ \ \ \ }&{T_j,}& {V_{lm},} &{ V_{km}, }&{ V_{kl}\ }.
}
$$
We shall express the coordinates of an element $y\in X_{\RT}\simeq \CZ^4$ as follows:
\bea(cc)
{(U_j,T_k,T_l,T_m)=({u_j},t_k,t_l,t_m),}&\hskip 35pt \RT={\Delta_j},\cr
{(U_j^{-1},V_{jk},V_{jl},V_{jm})=(\overline{u_j},v_{jk},v_{jl},v_{jm}),}&\hskip 35pt  \RT=C_j,\cr
{(T_j,V_{lm},V_{km},V_{kl})=(t_j,v_{lm},v_{km},v_{kl}),}&\hskip 35pt \RT=C^\prime_j,
\elea(eq:coord4)
and the corresponding ideal $I(y)$ is given by
\bea(cc)
{\langle Z_j-u_j Z_kZ_lZ_m, Z_k^2-t_k, Z_l^2-t_l, Z_m^2-t_m\rangle,} & {\RT=\Delta_j;} \cr 
\left.
\matrix{
\langle  Z_kZ_lZ_m-\overline{u_j}Z_j, Z_jZ_k-v_{jk}Z_lZ_m, Z_jZ_l-v_{jl}Z_kZ_m, \hfill\cr
 Z_jZ_m-v_{jm}Z_kZ_l, Z_j^2-t_j, Z_k^2-t_k, Z_l^2-t_l, Z_m^2-t_m \ \rangle,\ \hfill \cr
\left(t_j=\overline{u_j}v_{jk}v_{jl}v_{jm},\ t_\alpha=\overline{u_j}v_{j\alpha},\ (\alpha\ne j)\right),
}\right\}
&{\RT=C_j;} \cr
\left.
\matrix{
\langle Z_lZ_m-v_{lm}Z_jZ_k, Z_kZ_m-v_{km}Z_jZ_l,\hfill\cr 
Z_kZ_l-v_{kl}Z_jZ_m, Z_j^2-t_j, Z_k^2-t_k, Z_l^2-t_l, Z_m^2-t_m \rangle,\ \hfill\cr
\ (t_k=t_jv_{km}v_{kl},\ t_l=t_jv_{lm}v_{kl},\ t_m=t_jv_{km}v_{lm}),
}\right\}
&{\RT=C_j^\prime.}
\elea(eq:Iy)

In each case, one can show that the generators in the expression of (\ref{eq:Iy}) form the reduced Gr\"{o}bner basis of $I(y)$, and $\lt_\prec(I(y))=I(x_\RT)$ for a weight order $\prec$ with the weight in Interior$(\sigma_\RT)$. For $\RT={\Delta_j}$, we have ${I(x_{\Delta_j})=\langle Z_j,Z_k^2,Z_l^2,Z_m^2\rangle}$, and 
$$I(x_{\Delta_j})^\bot=\{1,Z_k,Z_l,Z_m,Z_lZ_m,Z_kZ_m,Z_kZ_l,Z_kZ_lZ_m\}$$
gives rise to a $G$-regular monomial basis for $\CZ[Z]/I(y)$. Similarly, for $\RT=C_j, C_j^\prime$, we have 
$I(x_{C_j})=\langle Z_kZ_lZ_m, Z_jZ_k, Z_jZ_l, Z_jZ_m, Z_1^2,Z_2^2, Z_3^2,Z_4^2\rangle$ and \\
$I(x_{C^\prime_j})=\langle Z_lZ_m, Z_kZ_m, Z_kZ_l,Z_1^2,Z_2^2,Z_3^2,Z_4^2\rangle$ 
with
$$I(x_{C_j})^\bot=\{1,Z_1,Z_2,Z_3,Z_4, Z_lZ_m, Z_kZ_m, Z_kZ_l \},$$ $$I(x_{C^\prime_j})^\bot=\{1,Z_1,Z_2,Z_3,Z_4, Z_jZ_k, Z_jZ_l, Z_jZ_m \},$$
which give rise to corresponding $G$-regular monomial basis of $\CZ[Z]/I(y)$. 
Since $X_{\Xi^*}$ is birational over $S_G$, with the vector bundle ${\cal F}_{X_{\Xi^*}}$, we have an epimorphism 
$$\lambda:X_{\Xi^*}\rightarrow \hl^{A_1(4)}(\CZ^4)$$
with $I(\lambda(y))=I(y)$ for $y\in X_{\Xi^*}$. We are going to show the injectivity of $\lambda$. For each $\RT\in\Xi^*(3)$, the coordinates of an element in $X_\RT$ can be represented in the form $(p_i/q_i=\gamma_i)_{i=1}^4$ for some monomials $p_i$, $q_i$ described in (\ref{eq:coord4}). Let $y$, $y^\prime$ be two points in $X_\RT$ with the coordinates $(p_i/q_i=\gamma_i)_{i=1}^4$ and $(p_i/q_i=\gamma_i^\prime)_{i=1}^4$ respectively.  By (\ref{eq:coord4}), (\ref{eq:Iy}), we have $p_i-\gamma_iq_i\in I(y)$ (resp. $p_i-\gamma_i^\prime q_i\in I(y^\prime)$) with $q_i\in I(x_\RT)^\bot$ for each $i$. When $I(y^\prime)=I(y)$, we have $(\gamma_i^\prime-\gamma_i)q_i\in I(y)$. By $q_i\in I(x_\RT)^\bot$, one has $\gamma_i^\prime=\gamma_i$ for all $i$, which implies $y=y^\prime$. Hence, $\lambda$ in injective on each $X_\RT$. Furthermore, the ideal $I(y)$ for $y\in X_\RT$ is completely determined by its toric coordinates, which appear as some generating elements of $I(y)$ in (\ref{eq:Iy}). By the construction of toric variety, one can conclude that if $y\in X_\RT$ and $y^\prime\in X_{\RT^\prime}$ with $I(y)=I(y^\prime)$, the toric coordinates of $y$  for $X_\RT$ and $y^\prime$ for $X_{\RT^\prime}$ are related by the transition function of these affine charts, hence $y=y^\prime$. For example, for $y\in X_{C_4^\prime}$ and $y^\prime\in X_{C_1}$, the affine coordinate of $y\in X_{C_4^\prime}$ is given by 
$$(\xi_1,\xi_2,\xi_3,\xi_4):=\left( Z_4^2,\ {{Z_2Z_3}\over{Z_1Z_4}},\   {{Z_1Z_3}\over{Z_2Z_4}},\ {{Z_1Z_2}\over{Z_3Z_4}} \right)$$
and $y^\prime \in X_{C_1}$ is 
$$(\eta_1,\eta_2,\eta_3,\eta_4):=\left({{Z_2Z_3Z_4}\over{Z_1}},\  {{Z_1Z_4}\over{Z_2Z_3}},\ {{Z_1Z_3}\over{Z_2Z_4}},\ {{Z_1Z_2}\over{Z_3Z_4}} \right).$$
By the expression of $\xi$'s and $\eta$'s in terms of $Z_i$'s, we obtain the transition function on $X_{C_4^\prime}\cap X_{C_1}$:
\begin{equation}
(\xi_1,\xi_2,\xi_3,\xi_4)=(\eta_1\eta_2,\ \eta_2^{-1},\eta_3,\eta_4).
\label{eq:trel}
\end{equation}
Therefore, $\lambda$ defines an isomorphism between $X_{\Xi^*}$ and $\hl^G(\CZ^4)$. 

By (\ref{cb}), the canonical bundle of $X_{\Xi^*}$ is given by
$\omega_{X_{\Xi^*}} = {\cal O}_{X_{\Xi^*}}(E)$
with $E=D_c$. As the star of $c$ in $\Xi^*$ is given by the octahedron in [Fig. 1], we have $E\simeq(\PZ^1)^3$. One can apply the toric technique to determine the $(-1)$-hyperplane structure of fiber $\PZ^1$ of the $(k,l)$-th factor projection, $p_j: E\rightarrow (\PZ^1)^2$ for the normal bundle of $E$. For example, in the case of the projection of $E$
onto $(\PZ^1)^2$ corresponding to the $2$-convex set $\langle v^{12}, v^{13}, v^{34}, v^{24}\rangle$, the relation (\ref{eq:trel}) between the local coordinates of $X_{C_1}$ and $X_{C_4^\prime}$ implies that the restriction of the normal bundle of $E$ on each fiber $\PZ^1$ over $(3,4)$-plane is the $(-1)$-hyperplane bundle. This proves Theorem \ref{th:A1(4)}.\\ {\hfill QED}

By the property of the normal bundle of the canonical divisor $E$ on $\PZ^1$-fibers in Theorem \ref{th:A1(4)}, one can blow down $E$ to obtain three different crepant resolutions of $S_G$. In fact, all these crepant resolutions are toric varieties described as follows:  Let $\Xi_j$ be the refinement
of $\Xi$ by adding the segment connecting
$v^{j4}$ and $v^{kl}$ to divide the central
polygon $\Diamond$ into $4$ simplices, where $\{j,k,l \}= \{1,2,3 \}$. Then, we have the relation of refinements: $\Xi \prec \Xi_j \prec \Xi^* $ for $j=1,2,3$, and each $X_{\Xi_j}$ is a crepant resolution of $X_\Xi$.
For the corresponding decompositions of the
central core $\Diamond$, the refinement relations  
are given by 
$\Diamond \prec \Diamond_j \prec \Diamond^*$, $(1\le j\le 3)$, with the pictorial realization [Fig. 3]. 
\begin{figure}[ht]
$$\mbox{\psfig{figure=crep2.ps}}$$
$$\mbox{\psfig{figure=octa.ps}}\hskip 45pt \mbox{\psfig{figure=crep3.ps}}\hskip 45pt \mbox{\psfig{figure=hilbs.ps}}$$
$$\mbox{\psfig{figure=crep1.ps}}$$
\centerline{[Fig. 3] Toric representation of flops of $4$-folds over the isolated singularity} \\
\centerline{and dominated by the smooth 4-fold.}
\end{figure}
By the relation between $X_{\Diamond^*}$ and $X_{\Diamond}$, one can conclude that $x_\Diamond$ is the isolated singularity of $X_\Diamond$, with the equation (\ref{4sing}). The three $4$-folds $X_{\Diamond_i}$ are all
``small''\footnote{Here a ``small'' resolution means a resolution with the exceptional locus of codimension $\geq 2$.} resolutions of the singular variety $X_\Diamond$. The relationship between these three crepant resolutions of the isolated singularity is defined to be the flop of $4$-folds. Hence, we have shown
the following result:
\begin{thm} \label{th:flip}
There are three crepant resolutions of
$S_{ A_1(4)}$ obtained by blowing down the canonical divisor $E$ of ${\rm
Hilb}^{ A_1(4)}(\CZ^4)$ in Theorem
$\ref{th:A1(4)}$. Any two such resolutions differ
by a flop of $4$-folds. 
\end{thm}
{\bf Remark}: The ${\goth S}_4$-action on coordinates leaves $\Diamond$ and $\Diamond^*$ invariant, but permutes three $\Diamond_i$'s. A similar phenomenon will also appear in the case $G=A_1(5)$ in the next section.

\end{section}

\begin{section}{$A_1(5)$-Hilbert Scheme and Crepant Resolutions of $\CZ^5/A_1(5)$}
In this section, we shall discuss the case of $G=A_1(5)$. Throughout this section, the indices $i$, $j$, $k$, $l$, $m$ always mean a permutation of $(1,2,3,4,5)$. In addition to the elements $e^i$, $v^{ij}$ of $N$ in $\Delta$, we consider the following rational points:
$$u^i={{1}\over{4}}\sum_{\alpha\ne i} e^\alpha \ , \ \ \ \ \ \ \ \  w^i={{1}\over{6}} \left(e^i+\sum_{\alpha=1}^5 e^\alpha\right).$$
Then, one has a refinement $\Xi^*$ of $\Xi$ with elements in $\Xi^*(4)$ given by:
$$
\Delta_i=\langle e^i,v^{ij},v^{ik},v^{il},v^{im} \rangle,\ \ \ {\bf I}_i=\langle w^i,v^{ij},v^{ik},v^{il},v^{im}  \rangle,\ \ \ {\bf II}_{ij}=\langle w^i,u^j,v^{ik},v^{il},v^{im} \rangle,$$
$${\bf III}_{jk}={\bf III}_{kj}=\langle u^j, u^k,v^{im},v^{lm},v^{il}  \rangle,\ \ \ {\bf IV}_{i,jk}={\bf IV}_{i,kj}=\langle  w^i,u^j,u^k,v^{il},v^{im}  \rangle,$$
$${\bf V}_{im}={\bf V}_{mi}=\langle  w^i,w^m,u^j,u^k,u^l,v^{im} \rangle, \ \ {\bf VI}:=\langle  u^\alpha,w^\alpha \rangle_{\alpha=1}^5.$$
There are 5, 5, 20, 10, 30, 10, 1 elements of $\Delta_i$, ${\bf I}_i$, ${\bf II}_{ij}$, ${\bf III}_{jk}$, ${\bf IV}_{i,jk}$, ${\bf V}_{im}$, ${\bf VI}$ respectively. The facet relations for these $4$-polytopes can be depicted graphically in the following diagram:
$$
\matrix{\hfill\Delta_i- {\bf I}_i - {\bf II}_{ij}&-{{\bf IV}_{i,jk}}- & {\bf V}_{im} - {\bf VI}\hfill\cr
{}&{|}&{}\cr
{}&{\bf III}_{jk}&{}
}$$
where a line describes the intersection of two polytopes with a common facet. Except the heptahedrons ${\bf V}_{im}$ and the decahedron {\bf VI}, the rest of the $4$-polytopes are all simplicial. 
\begin{thm} For $G=A_1(5)$,
$\hl^G(\CZ^5)$ is the toric variety $X_{\Xi^*}$ with the canonical sheaf
$\omega_{X_{\Xi^*}} = {\cal O}_{X_{\Xi^*}}
\left(\sum_{i=1}^5 (2D_{w^i}+D_{u^i}) \right)$, 
where $D_{u^i}$ and $D_{w^i}$ are the divisors corresponding to $u^i$, $w^i$ respectively. 
\label{Thm:Main}
\end{thm}
\noindent{\it Proof}: With (\ref{eq:TUV}) for $n=5$, we have
$U_i={{Z_i}\over {Z_jZ_kZ_lZ_m}}$, $V_{ij}={{Z_iZ_j}\over{Z_kZ_lZ_m}}$.  
The generators of $N \cap \sigma_{\RT}$ and $M \cap \check{\sigma_{\RT}}$ for $\RT\in \Xi^*(4)$ are given by the following table:
\vskip 10pt
\centerline{
\begin{tabular}{|c|c|c|}\hline
$\RT$ & $N \cap \sigma_{\RT}$ & $M \cap \check{\sigma_{\RT}}$ \\ \hline
$\Delta_i$ & $e^i,v^{ij},v^{ik},v^{il},v^{im}$  & $U_i,T_j,T_k,T_l,T_m$\\
${\bf I}_i$ & $3w^i,v^{ij},v^{ik},v^{il},v^{im}$ & $U^{-1}_i,V_{ij},V_{ik},V_{il},V_{im}$ \\
 ${\bf{II}}_{ij}$ & $3w^i,2u^j,v^{ik},v^{il},v^{im}$ & $V^{-1}_{ij},T_j,V_{ik},V_{il},V_{im}$ \\
 ${\bf{III}}_{jk}$ & $2u^j, 2u^k,v^{im},v^{lm},v^{il}$ & $T_j, T_k, V_{im}, V_{lm},  V_{il}$ \\
 ${\bf{IV}}_{i,jk}$ & $3w^i,2u^j,2u^k,v^{il},v^{im}$ & $V^{-1}_{lm}, V^{-1}_{ij}, V^{-1}_{ik}, V_{il}, V_{im}$ \\
 ${\bf V}_{im}$ & $3w^i,3w^m,2u^j,2u^k,2u^l,v^{im}$ & $V^{-1}_{ij}, V^{-1}_{jm},  V^{-1}_{ik}, V^{-1}_{km},V^{-1}_{il}, V^{-1}_{lm},V_{im}$ \\
 ${\bf VI}$ & $3w^i$, $2u^i$, $(1\le i\le 5)$ & $V^{-1}_{ij}$, $(1\le i<j\le 5)$\\ 
\hline
\end{tabular}
}
\vskip 10pt
\noindent Then, one can easily see that $X_\RT$ has the smooth toric structure $\CZ^5$ except $\RT={{\bf V}_{im}}$ and ${\bf VI}$, which contribute the singularities of $X_{\Xi^*}$. In fact, $X_{{\bf V}_{im}}$ is the affine variety in $\CZ^7$ defined by the relations,
\begin{equation}
V^{-1}_{ij}V^{-1}_{km}=V^{-1}_{ik}V^{-1}_{jm}\ ,\ \ \ 
V^{-1}_{ik}V^{-1}_{lm}=V^{-1}_{il}V^{-1}_{km}\ ,\ \ \ 
V^{-1}_{il}V^{-1}_{jm}=V^{-1}_{ij}V^{-1}_{lm}\ ;
\label{eq:relV}
\end{equation}
and $X_{\bf VI}$ is the affine variety in $\CZ^{10}$ given by
\begin{equation}
V^{-1}_{ij}V^{-1}_{il}V^{-1}_{km}= V^{-1}_{ik} V^{-1}_{im}V^{-1}_{jl},\  V^{-1}_{ij}V^{-1}_{kl}=V^{-1}_{ik}V^{-1}_{jl}
\label{eq:relVI}
\end{equation}
where the indices run through all possible $i,j,k,l,m$. For $\gamma\in \CZ$, we define the following eigen-polynomials of $G$:
$$F_{ij}(\gamma):=Z_iZ_j-\gamma Z_kZ_lZ_m, \ \ H_{ij}(\gamma):=Z_kZ_lZ_m-\gamma Z_iZ_j.$$ 

We shall denote the coordinates of an element $y\in X_\RT$ by
\vskip 6pt
\centerline{
\begin{tabular}{|c|c|}\hline
$\RT$ & $y\in X_\RT$ \\ \hline
$\Delta_i$ & 
$\left(U_i,T_j,T_k,T_l,T_m\right) =\left(\gamma_i,t_j,t_k,t_l,t_m\right)$\\
${\bf I}_i$ & $\left(U_i^{-1},V_{ij},V_{ik},V_{il},V_{im}\right) =\left(\gamma_{jklm}, \gamma_{ij},\gamma_{ik},\gamma_{il},\gamma_{im} \right)$ \\
 ${\bf{II}}_{ij}$ & 
$\left(V^{-1}_{ij},T_j,V_{ik},V_{il},V_{im}\right) =\left(\gamma_{klm}, t_j, \gamma_{ik}, \gamma_{il}, \gamma_{im}\right)$ \\
 ${\bf{III}}_{jk}$ & 
$\left( V_{im},V_{il},V_{lm}, T_k,T_j\right) =\left(\gamma_{im}, \gamma_{il}, \gamma_{lm}, t_k,t_j\right)$ \\
 ${\bf{IV}}_{i,jk}$ & 
$\left( V_{il}, V_{im}, V^{-1}_{ij}, V^{-1}_{ik}, V^{-1}_{lm}\right) =\left(\gamma_{il}, \gamma_{im}, \gamma_{klm}, \gamma_{jlm}, \gamma_{ijk}\right)$ \\
 ${\bf V}_{im}$ & 
$\matrix{
 &\left(V_{im},V^{-1}_{il}, V^{-1}_{lm}, V^{-1}_{ik} , V^{-1}_{ij}, V^{-1}_{jm}, V^{-1}_{km}\right)\\
{=}&\left(\gamma_{im}, \gamma_{jkm}, \gamma_{ijk}, \gamma_{jlm}, \gamma_{klm}, \gamma_{ikl}, \gamma_{ijl}\right)
}$
\\
 ${\bf VI}$ & $\left( V^{-1}_{ij}\right)_{i<j}= \left(\gamma_{klm}\right)_{k<l<m}$ \\ 
\hline
\end{tabular}
}
\vskip 6pt

\noindent where the $\gamma$'s and $t$'s are complex numbers, and in the cases of ${\bf V}_{im}$ and ${\bf VI}$, the coordinates are governed by the relations (\ref{eq:relV}), (\ref{eq:relVI}) respectively. Using the above coordinates of $y$, one can express the ideal $I(y)$ as follows:

\begin{tabular}{cc}
$\Delta_i$ & 
$I(y)=\langle Z_j^2-t_j, Z_k^2-t_k, Z_l^2-t_l, Z_m^2-t_m, Z_i-\gamma_i Z_j Z_k Z_l Z_m  \rangle$,
\\
\\
${\bf I}_i$ 
&$\left\{
\matrix{
  \matrix{
    I(y) & = &\langle 
           Z_\alpha^2-t_\alpha
         \rangle_{\alpha=1}^5+ \langle F_{ij}(\gamma_{ij}),
            F_{ik}(\gamma_{ik}), \hfill\\
       {}&{}&F_{il}(\gamma_{il}), 
             F_{im}(\gamma_{im}),
             Z_jZ_kZ_lZ_m-\gamma_{jklm}Z_i 
         \rangle\hfill
  }\\
  \left(\matrix{
          t_j=\gamma_{ij}\gamma_{jklm},\ 
          t_k=\gamma_{ik}\gamma_{jklm},\ 
          t_l=\gamma_{il}\gamma_{jklm},\\ 
          t_m=\gamma_{im}\gamma_{jklm},\ 
          t_i=\gamma_{ij}\gamma_{ik}\gamma_{il}\gamma_{im}\gamma_{jklm}^2 
        }
  \right)
}
\right.$
\\
\\
${\bf II}_{ij}$ & 
$\left\{
\matrix{
  {I(y)}=\langle 
         Z_\alpha^2-t_\alpha
       \rangle_{\alpha=1}^5+
       \langle
         F_{ik}(\gamma_{ik}), F_{il}(\gamma_{il}), F_{im}(\gamma_{im}), 
         H_{ij}(\gamma_{klm})
       \rangle  \\
  \left(\matrix{
         t_k=\gamma_{ik}\gamma_{klm}t_j,\ 
         t_l=\gamma_{il}\gamma_{klm}t_j,\ 
         t_m=\gamma_{im}\gamma_{klm}t_j,\\  
         t_i=\gamma_{ik}\gamma_{il}\gamma_{im}{t_j}^2\gamma_{klm}
        }
  \right)
}
\right.$
\\
\\
${\bf III}_{jk}$ & 
$\left\{
\matrix{
  I(y)=\langle 
         Z_\alpha^2-t_\alpha
       \rangle_{\alpha=1}^5+
       \langle 
         F_{im}(\gamma_{im}), F_{il}(\gamma_{il}), F_{lm}(\gamma_{lm})
       \rangle, 
  \\
  \left(t_i=\gamma_{im}\gamma_{il}t_kt_j,\ 
        t_l=\gamma_{il}\gamma_{lm}t_kt_j,\ 
        t_m=\gamma_{im}\gamma_{lm}t_kt_j\ 
  \right)
}
\right.$
\\
\\
${\bf IV}_{i,jk}$ &
$\left\{
\matrix{
 \matrix{
  {I(y)}&{=}&
       \langle 
         Z_\alpha^2-t_\alpha
       \rangle_{\alpha=1}^5+ 
       \langle
         F_{im}(\gamma_{im}), F_{il}(\gamma_{il}), 
       \hfill \\
      {}&{}&
          H_{ij}(\gamma_{klm}), H_{ik}(\gamma_{jlm}), H_{lm}(\gamma_{ijk}) 
      \rangle \hfill
 }
\\
  \left(
    \matrix{
      t_m= \gamma_{im} \gamma_{jlm}\gamma_{klm} \gamma_{ijk},\ 
      t_l= \gamma_{il} \gamma_{jlm}\gamma_{klm} \gamma_{ijk},\ 
      \\
      t_i=\gamma_{im} \gamma_{il} \gamma_{jlm}\gamma_{klm} {\gamma_{ijk}}^2,
      t_j= \gamma_{jlm}\gamma_{ijk},\ t_k= \gamma_{lmk}\gamma_{ijk}
    }
  \right)
}
\right.$
\\
\\
${\bf V}_{im}$ &
$\left\{
\matrix{
 \matrix{
  {I(y)}&{=}&
       \langle 
         Z_\alpha^2-t_\alpha\rangle_{\alpha=1}^5
         +\langle F_{im}(\gamma_{im}),
                  H_{il}(\gamma_{jkm}), H_{ik}(\gamma_{jlm}), \hfill  \\
      {}&{}&
         H_{ij}(\gamma_{klm}), H_{jm}(\gamma_{ikl}), H_{lm}(\gamma_{ijk}),
         H_{km}(\gamma_{ijl}) 
       \rangle\hfill
 }
\\
  \left(
    \matrix{
      t_i=\gamma_{{im}}\gamma_{{ikl}}\gamma_{{jlm}}\gamma_{{ijk}}
         =\gamma_{{im}}\gamma_{{ikl}}\gamma_{{jkm}}\gamma_{{ijl}},\ \\
      t_m=\gamma_{im}\gamma_{{klm}}\gamma_{{jlm}}\gamma_{{ijk}}
         =\gamma_{{im}}\gamma_{{klm}}\gamma_{{jkm}}\gamma_{{ijl}},\ \\
      t_l=\gamma_{{ikl}}\gamma_{{jlm}}=\gamma_{ijl}\gamma_{klm},\ 
      t_j=\gamma_{{ijk}}\gamma_{{jlm}}=\gamma_{ijl}\gamma_{jkm},\ \\
      t_k=\gamma_{{ijk}}\gamma_{{klm}}=\gamma_{ikl}\gamma_{jkm}\ 
    }
  \right)
}
\right.$
\\
\\
${\bf VI}$ &
$\left\{
\matrix{
  I(y)= \langle 
          Z_\alpha^2-t_\alpha 
        \rangle_{\alpha=1}^5+
        \langle H_{lm}(\gamma_{ijk})
          \vert {\rm \ for\ all \ possible\ \ } i,j,k,l,m 
        \rangle,\\
  \left(\matrix{
         t_i=\gamma_{ijk}\gamma_{ilm}=\gamma_{ijm}\gamma_{ilk}, 
        \gamma_{klm}\gamma_{jkm}\gamma_{ijl}=
        \gamma_{jlm}\gamma_{jkl}\gamma_{ikm} \\
        {\rm \ for\ all \ possible\ \ } i,j,k,l,m
        }
  \right)
}
\right.$
\\
\end{tabular}

\noindent One can show that the generators in the above expressions of $I(y)$ form the corresponding reduced Gr\"{o}bner basis with respect to a weight order $\prec$ with weight ${\bf w}\in{\rm Interior}(\sigma_\RT)$. Furthermore, the ideal $\lt_\prec(I(y))$ is equal to $I(x_\RT)$. Therefore, $I(x_\RT)^\bot$ gives rise to a basis of $\CZ[Z]/I(y)$, which is $G$-regular by the following explicit description of basis elements. 
\vskip 6pt
\centerline{
\begin{tabular}{|c|l|}\hline
$\RT$ &\hskip 100pt ${I(x_\RT)}^\bot$ \\ \hline
$\Delta_i$ &
$\matrix{
{1,Z_j,Z_k,Z_l,Z_m,Z_jZ_k,Z_jZ_l,Z_kZ_l,Z_kZ_m,Z_lZ_m,Z_jZ_m,}\\
{Z_jZ_kZ_l,Z_jZ_kZ_m, Z_jZ_lZ_m,Z_kZ_lZ_m,Z_jZ_kZ_lZ_m   }\hfill
}$ \\ \hline
${\bf I}_i $ & $\matrix{
{1,Z_i,Z_j,Z_k,Z_l,Z_m,Z_jZ_k,Z_jZ_l,Z_kZ_l,Z_kZ_m,Z_lZ_m,Z_jZ_m, }\\{Z_jZ_kZ_l,Z_jZ_kZ_m,Z_jZ_lZ_m,Z_kZ_lZ_m}\hfill }$\\ \hline
${\bf II}_{ij}$ &  $\matrix{
{1,Z_i,Z_j,Z_k,Z_l,Z_m,Z_iZ_j,Z_jZ_k,Z_jZ_l,Z_kZ_l,Z_kZ_m,Z_lZ_m,}\\
 Z_jZ_m,Z_jZ_kZ_l,Z_jZ_kZ_m,Z_jZ_lZ_m\hfill 
}
$\\ \hline
${\bf III}_{jk}$ & $\matrix{
1,Z_i,Z_j,Z_k,Z_l,Z_m,Z_iZ_j,Z_jZ_k,Z_jZ_l,Z_kZ_l,Z_kZ_m,Z_iZ_jZ_k,\\ 
Z_jZ_m, Z_jZ_kZ_l,Z_jZ_kZ_m,Z_iZ_k\hfill}$\\ \hline
${\bf IV}_{i,jk}$ & $\matrix{
1,Z_i,Z_j,Z_k,Z_l,Z_m,Z_iZ_j,Z_jZ_k,Z_jZ_l,Z_kZ_l,Z_kZ_m,Z_lZ_m,\\ 
Z_jZ_m,Z_jZ_kZ_l,Z_jZ_kZ_m,Z_iZ_k\hfill} $\\ \hline
${\bf V}_{im}$ & $\matrix{ 1,Z_i,Z_j,Z_k,Z_l,Z_m,Z_iZ_j,Z_jZ_k,Z_jZ_l,Z_kZ_l,Z_kZ_m,Z_iZ_jZ_k,\\
Z_jZ_m, Z_jZ_kZ_l, Z_iZ_l,Z_iZ_k \hfill}
$\\ \hline
${\bf VI}$ & $\matrix{ 1,Z_1,Z_2,Z_3,Z_4,Z_5,Z_1Z_2,Z_1Z_3,Z_1Z_4,Z_1Z_5,Z_2Z_3,Z_2Z_4,\\
Z_2Z_5,Z_3Z_4,Z_3Z_5,Z_4Z_5\hfill}$\\
\hline
\end{tabular}}
\vskip 6pt
\noindent Thus, $\CZ[Z]/I(y)$ is a regular $G$-module for each $y\in X_{\Xi^*}$, and $X_{\Xi^*}$ is birational over $S_G$ having the vector bundle ${\cal F}_{\Xi^*}$. Hence, there is an epimorphism $\lambda:X_{\Xi^*}\rightarrow \hl^G(\CZ^5)$ with $I(\lambda(y))=I(y)$. With the same argument as in the case $n=4$, one can show that $\lambda$ is injective on each $X_\RT$ by the reason that the toric coordinates are encoded in the generators of $I(y)$. By the toric variety structure of $X_{\Xi^*}$, one obtains the injectivity of $\lambda$. Hence $X_{\Xi^*}\simeq \hl^G(\CZ^5)$.
By (\ref{cb}), one has the expression of the canonical sheaf for $X_{\Xi^*}$ .\\
{\hfill QED.}

We are going to construct certain crepant resolutions of $S_G$ dominated by the toric variety $X_{\Xi^*}\simeq\hl^G(\CZ^5)$. Consider the permutation $\tau=(12345)$ on the coordinates of $N_0$. Define the following simplices in $\Delta$:
$$C:=\langle v^{12},v^{23},v^{34},v^{45},v^{15}\rangle,$$
$$D_0:=\langle v^{12},v^{23},v^{34},v^{13},v^{15}\rangle,$$
$$E_0:=\langle v^{35},v^{23},v^{25},v^{45},v^{15}\rangle,$$
and $D_i:=\tau^i(D_0)$, $E_i:=\tau^i(E_0)$ ($i=1,..,4$). Then, one can show that $C$, $D_i$ and $E_i$ ($0\le i\le 4$) form an integral simplicial decomposition of the central core $\Diamond\in \Xi(4)$, hence one has the refinement $\Diamond^\prime$ of $\Diamond$ with $\Diamond^\prime(4)=\{C\} \cup\{D_i,E_i\}_{i=0}^4$. Denote $\Xi^\prime$ the rational polytope decomposition of $\Delta$ with $\Xi^\prime(4)=\Diamond^\prime(4)\cup\{\Delta_i\}_{i=1}^5$. Then $\Xi^\prime$ is a refinement of $\Xi$. By computing the primitive generators of each element in $\Xi^\prime(4)$, one can easily check that $X_{\Xi^\prime}$ is a crepant resolution of $S_G$. The facet relations of simplices in $\Xi^\prime(4)$ are given in the following diagram: 
$$
{\put(13,40){\line(1,0){12}}
\put(-13,40){\line(-1,0){12}}
\put(33,35){\line(1,-3){4}}
\put(-33,35){\line(-1,-3){4}}
\put(44,6){\line(1,-3){4}}
\put(-44,6){\line(-1,-3){4}}
\put(-34,-30){\line(-4,3){10}}
\put(34,-30){\line(4,3){10}}
\put(10,-48){\line(4,3){12}}
\put(-10,-48){\line(-4,3){12}}
\put(0,12){\line(0,1){20}}
\put(9,3){\line(3,1){19}}
\put(9,-10){\line(3,-4){12}}
\put(-9,3){\line(-3,1){19}}
\put(-9,-10){\line(-3,-4){12}}
\put(0,-65){\line(0,-1){12}}
\put(63,-22){\line(3,-1){10}}
\put(-63,-22){\line(-3,-1){10}}
\put(36,48){\line(3,4){10}}
\put(-36,48){\line(-3,4){10}}
\put(45,63){$\Delta_4$}
\put(-55,63){$\Delta_1$}
\put(-85,-32){$\Delta_3$}
\put(73,-32){$\Delta_2$}
\put(-4,-4){$C$}
\put(-8,36){$D_0$}
\put(-50,9){$D_2$}
\put(34,9){$D_3$}
\put(27,37){$E_4$}
\put(-40,37){$E_1$}
\put(-10,-60){$E_{0}$}
\put(-5,-85){$\Delta_5$}
\put(48,-20){$E_2$}
\put(-60,-20){$E_3$}
\put(22,-40){$D_1$}
\put(-36,-40){$D_4$}}
$$
\begin{lma} We have the relation
$\Diamond\prec\Diamond^\prime\prec\Diamond^*$, consequently, $\Xi\prec\Xi^\prime\prec\Xi^*$.
\label{lma:refine}
\end{lma}
\noindent{\it Proof}: It suffices to show the relation $\Diamond^\prime\prec\Diamond^*$. By a detailed analysis of the $4$-polytopes in $\Xi^*(4)$, one can conclude that the simplex $C$ is the union of the following eleven polytopes in $\Xi^\prime(4)$ with the facet relations: 
$$
\put(-100,-4){$C:$}
\put(13,40){\line(1,0){12}}
\put(-13,40){\line(-1,0){12}}
\put(33,35){\line(1,-3){4}}
\put(-33,35){\line(-1,-3){4}}
\put(44,6){\line(1,-3){4}}
\put(-44,6){\line(-1,-3){4}}
\put(-34,-30){\line(-4,3){10}}
\put(34,-30){\line(4,3){10}}
\put(10,-48){\line(4,3){12}}
\put(-10,-48){\line(-4,3){12}}
\put(0,12){\line(0,1){20}}
\put(9,3){\line(3,1){19}}
\put(9,-10){\line(3,-4){12}}
\put(-9,3){\line(-3,1){19}}
\put(-9,-10){\line(-3,-4){12}}
\put(-8,-4){${\bf{VI}}$}
\put(-8,36){${{\bf{V}}}_{34}$}
\put(-50,9){${{\bf{V}}}_{23}$}
\put(34,9){${{\bf{V}}}_{45}$}
\put(31,37){${{\bf{IV}}}_{412}$}
\put(-54,37){${{\bf{IV}}}_{315}$}
\put(-10,-60){${{\bf{IV}}}_{134}$}
\put(48,-20){${{\bf{IV}}}_{523}$}
\put(-75,-20){${{\bf{IV}}}_{245}$}
\put(22,-40){${{\bf{V}}}_{15}$}
\put(-40,-40){${{\bf{V}}}_{12}$}
$$
Similarly, the relations of $D_0$, $E_0$ with $4$-polytopes in $\Diamond^*(4)$ are given in the following diagram:
$$
\put(-160,-4){$D_0:$} 
\put(-70,9){\line(0,1){15}}
\put(-70,-9){\line(0,-1){15}}
\put(-57,0){\line(1,0){27}}
\put(-83,0){\line(-1,0){13}}
\put(-83,30){\line(-1,0){13}}
\put(-83,-30){\line(-1,0){13}}
\put(-44,-30){\line(1,0){14}}
\put(-120,9){\line(0,1){15}}
\put(-120,-9){\line(0,-1){15}}
\put(-20,-9){\line(0,-1){15}}
\put(-77,-4){${\bf{V}}_{13}$}
\put(-127,-4){${\bf{IV}}_{145}$}
\put(-27,-4){${\bf{IV}}_{325}$}
\put(-77,-34){${\bf{IV}}_{345}$}
\put(-77,26){${\bf{IV}}_{124}$}
\put(-127,26){${\bf{II}}_{14}$}
\put(-127,-34){${\bf{III}}_{45}$}
\put(-27,-34){${\bf{II}}_{35}$}
\put(20,-4){$E_0:$}
\put(100,5){\line(0,-1){15}}
\put(120,14){\line(1,0){15}}
\put(87,14){\line(-1,0){13}}
\put(87,-16){\line(-1,0){13}}
\put(60,5){\line(0,-1){15}}
\put(93,10){${\bf{IV}}_{514}$}
\put(53,10){${\bf{II}}_{51}$}
\put(143,10){${\bf{III}}_{14}$}
\put(93,-20){${\bf{II}}_{54}$}
\put(58,-20){${\bf I}_{5}$}
$$
By applying $\tau^i$ ($1\le i\le 4$) on the above diagrams for $D_0$ and $E_0$, one obtains the decompositions of $D_i$, $E_i$ in terms of elements in $\Diamond^*(4)$. The results then follow. \\
{\hfill QED} 

As known in \S2, $\Xi$ and $\Xi^*$ are ${\goth S}_5$-invariant polytope decompositions of $\Delta$, but not for $\Xi^*$. The isotropy subgroup of ${\goth S}_5$ for $\Xi^\prime$ is ${\goth G}^\prime=\langle (25)(34),(12345)\rangle$, which is isomorphic to the dihedral group of order $10$, with the index $[{\goth S}_5:{\goth G}^\prime]=12$. By applying permutation elements in ${\goth S}_5$ to the decomposition $\Xi^\prime$, one obtains twelve integral simplicial refinements $\Xi_i$ ($1\le i\le 12$) of $\Xi$, with $\Xi\prec\Xi_i\prec\Xi^*$. Correspondingly, there are twelve decompositions $\Diamond_i$ ($1\le i\le 12$) of the central core $\Diamond$ with the following refinement relations: 
$$\Diamond\prec\Diamond_i \prec \Diamond^*\ ,\ \ \ \ i=1,..,12.$$
The connection between these twelve smooth $5$-folds corresponding to the toric data $\Diamond_i$'s can be regarded as the ``flop'' of $5$-folds, all of which are crepant resolutions of the singular variety defined by (\ref{5sing}). Therefore, we have obtained the following result:
\begin{thm} There are twelve toric crepant resolutions, $X_{\Xi_i}$ $1\le i\le 12$, of $S_{A_1(5)}$  which are dominated by $\HSA$. Any two such resolutions differ by a ``flop'' of $5$-folds.
\label{Thm:Flops}
\end{thm}

\end{section}

\end{document}